\documentclass[11pt]{article}
\usepackage{mathrsfs,amsmath,amssymb,amsfonts}
\numberwithin{equation}{section}
\begin{document}
\title{Another weighted approximation of functions with singularities by combinations of Bernstein operators}
\author{Wen-ming Lu$^1$, Lin Zhang$^2$\thanks{Corresponding author. Address: Department of Mathematics, Zhejiang University, Hangzhou 310027, PR China. \textsl{E-mail address}:
\texttt{linyz@zju.edu.cn}(L.-Zhang);\texttt{lu\_wenming@163.com}(W.-Lu).}}
\date{\it $^1$School of Science, Hangzhou Dianzi University,
Hangzhou, 310018 P.R. China\\$^2$Department of Mathematics, Zhejiang
University, Hangzhou, 310027 P.R. China} \maketitle
\mbox{}\hrule\mbox{}\\[0.5cm]
\textbf{Abstract}\\[-0.2cm]^^L
A new type of combinations of Bernstein operators is given in
\cite{S. Yu}. Here, we introduce another one, which can be used to
approximate the functions with singularities.
The direct and inverse results of the weighted approximation of this new type combinations are given.
\mbox{}\hrule\mbox{}
\section{Introduction}
The present work continues to study modified Bernstein operators
following \cite{S. Yu}. Here, the notations are referred to \cite{S.
Yu}. For convenience, these notations will be listed. The set of all
continuous functions, defined on the interval $I$, is denoted by
$C(I)$. For any $f\in C([0,1])$, the corresponding Bernstein
operators are defined as follows:
$$B_n(f,x):=\sum_{k=0}^nf(\frac{k}{n})p_{nk}(x),$$
where
$$p_{nk}(x):={n \choose k}x^k(1-x)^{n-k}, \ k=0,1,2,\ldots,n, \ x\in[0,1].$$
Approximation properties of Bernstein operators have been studied
very well (see \cite{Berens}, \cite{Della},
\cite{Totik}-\cite{Lorentz}, \cite{Yu}-\cite{X. Zhou}, for example).
In order to approximate the functions with singularities, Della
Vecchia et al. \cite{Della} and Yu-Zhao \cite{Yu} introduced some
kinds of modified Bernstein
operators.\\~\\
Let $$\bar{w}(x)=|x-\xi|^\alpha,\ 0<\xi<1,\ \alpha>0,$$ and
$$C_{\bar{w}}:= \{{f \in C([0,1] \setminus {\xi})
:\lim\limits_{x\longrightarrow\xi}(\bar{w}f)(x)=0 }\}.$$ The norm in
$C_{\bar{w}}$ is defined as
$\|f\|_{C_{\bar{w}}}:=\|{\bar{w}}f\|=\sup\limits_{0\leqslant
x\leqslant 1}|({\bar{w}f})(x)|$. Define
$$W_{\bar{w}}^{r}:= \{f \in
C_{\bar{w}}:f^{(r-1)} \in A.C.((0,1)),\
\|{\bar{w}}\varphi^rf^{(r)}\|<\infty\}.$$ For $f \in C_{\bar{w}}$,
define the weighted modulus of smoothness
by\\
$$\omega_{\varphi}^{r}(f,t)_{\bar{w}}:=\sup_{0<h\leqslant
t}\{\|{\bar{w}}\triangle_{h\varphi}^{r}f\|_{[16h^2,1-16h^2]}+\|{\bar{w}}{\overrightarrow{\triangle}_{h}^{r}f\|_{[0,16h^2]}+\|{\bar{w}}{\overleftarrow{\triangle}_{h}^{r}}f\|_{[1-16h^2,1]}}\},$$
where
$$\Delta_{h\varphi}^{r}f(x)=\sum_{k=0}^{r}(-1)^{k}{r \choose k}f(x+(\frac
r2-k)h\varphi(x)),$$
$$\overrightarrow{\Delta}_{h}^{r}f(x)=\sum_{k=0}^{r}(-1)^{k}{r \choose k}f(x+(r-k)h),$$
$$\overleftarrow{\Delta}_{h}^{r}f(x)=\sum_{k=0}^{r}(-1)^{k}{r \choose k}f(x-kh),$$
and $\varphi(x)=\sqrt{x(1-x)}$. The weighted $K$-function is given
by
$$K_{r,\varphi}(f,t^r)_{\bar{w}}:=\underset{g}{\inf
}\{\|{\bar{w}}(f-g)\|+t^{r}\|{\bar{w}}\varphi^{r}g^{(r)}\|:g \in
W_{\bar{w}}^{r}\}.$$ It was shown in \cite{Totik} that
$K_{\varphi}(f,t^r)_{\bar{w}} \thicksim
\omega_{\varphi}^{r}(f,t)_{\bar{w}}$. Della Vecchia et al. firstly
introduced $B_{n}^{\ast}(f,x)$ and ${\bar{B}}_{n}(f,x)$ in
\cite{Della}, where the properties of $B_{n}^{\ast}(f,x)$ and
${\bar{B}}_{n}(f,x)$ are studied. Among others, they prove that
$$\|w(f-B_{n}^{\ast}(f))\|\leqslant C\omega_{\varphi}^{2}(f,n^{-1/2}),\ f\in
C_{w},$$
$$\|{\bar{w}}(f-{\bar{B}_{n}(f)})\|\leqslant
\frac{C}{n^{3/2}}\sum_{k=1}^{[\sqrt{n}]}k^{2}\omega_{\varphi}^{2}(f,\frac{1}{k})_{\bar{w}}^{\ast},\
f\in C_{\bar{w}},$$ where $w(x)=x^{\alpha}(1-x)^{\beta},\
\alpha,\ \beta\geqslant 0,\ \alpha+\beta>0,\ 0\leqslant x \leqslant1.$
In \cite{S. Yu}, for any $\alpha,\ \beta>0,\ n\geqslant
2r+\alpha+\beta$, there hold
$$\|wB_{n,r}^{\ast}(f)\|\leqslant C\|wf\|,\ f\in C_{w},$$
\begin{eqnarray*}
\|w(B_{n,r}^{\ast}(f)-f)\|\leqslant \left\{
\begin{array}{lrr}
{\frac{C}{n^{r}}} (\|wf\|+\|w\varphi^{2r}f^{(2r)}\|),  &&f\in
W_{w}^{2r},    \\
C(\omega_{\varphi}^{2r}(f,n^{-1/2})_{w}+n^{-r}\|wf\|),   &&f\in
C_{w},
              \end{array}
\right.
\end{eqnarray*}
\begin{eqnarray*}
\|w\varphi^{2r}B_{n,r}^{\ast(2r)}(f)\|\leqslant \left\{
\begin{array}{lrr}
Cn^{r}\|wf\|,    &&f\in C_{w},    \\
C(\|wf\|+\|w\varphi^{2r}f^{(2r)}\|),    &&f\in W_{w}^{2r}.
              \end{array}
\right.
\end{eqnarray*}
and for $0< \gamma <2r,$
$$\|w(B_{n,r}^{\ast}(f)-f)\|=O(n^{-\gamma/2}) \Longleftrightarrow
\omega_{\varphi}^{2r}(f,t)_{w}=O(t^{r}).$$ On the other hand, since
the Bernstein polynomials cannot be used for the investigation of
higher orders of smoothness, Butzer \cite{Butzer} introduced the
combinations of Bernstein polynomials which have higher orders of
approximation. Ditzian and Totik \cite{Totik} extended this method
of combinations and defined the following combinations of Bernstein
operators:
\begin{eqnarray}
B_{n,r}(f,x):=\sum_{i=0}^{r-1}C_{i}(n)B_{n_i}(f,x)\label{s1}
\end{eqnarray}
with the
conditions
\\$(a)n=n_0<n_1< \cdots <n_{r-1}\leqslant Cn,$
\\$(b)\sum_{i=0}^{r-1}|C_{i}(n)|\leqslant C,$
\\$(c)\sum_{i=0}^{r-1}C_{i}(n)=1,$
\\$(d)\sum_{i=0}^{r-1}C_{i}(n)n_{i}^{-k}=0$,\ for $k=1,\ldots,r-1$.
\\Some approximation behaviors of the operators defined as (\ref{s1}) can
be found in \cite{Ditzian}-\cite{Guo}, \cite{Xie} and \cite{L.S.
Xie}. For example, Ditzian and Totik \cite{Totik} showed that
$$\|B_{n,r}(f)-f\|\leqslant C(\omega_{\varphi}^{2r}(f,n^{-1/2})+n^{-r}\|f\|),$$
and for $0< \alpha <2r,$
$$\|B_{n,r}(f)-f\|=O(n^{-\alpha /2}) \Longleftrightarrow
\omega_{\varphi}^{2r}(f,n^{-1/2})=O(t^{\alpha}),$$ where
$\omega_{\varphi}^{2r}(f,t)$ is the modulus of smoothness with the
step-weight function $\varphi(x),$ and $\|f\|=\|f\|_{C([0,1])}$.
\\~\\
The main purpose of the present paper is to give another new type of
combinations of Bernstein operators (combinations defined as
(\ref{s1}) cannot be used to approximate functions in $C_{\bar{w}}$)
so as to obtain higher approximation order. In Section 2, we will
give the new type of combinations, and the direct and inverse
results of the weighted approximation by the new type of
combinations. Some lemmas will be given in Section 3, while the
proofs of the results will be given in Section 4. Throughout the
paper, $C$ denotes a positive constant independent of $n$ and $x$,
which may be different in different cases.
\section{The main results}
For any positive integer $r$, we consider the determinant
\begin{eqnarray*}
A_{r}:=
\begin {matrix}
\begin{vmatrix}
1 & 1 & 1 & \cdots & 1 \\
2r+1 & 2r+2 & 2r+3 & \cdots & 4r+1 \\
(2r)(2r+1) & (2r+1)(2r+2) & (2r+2)(2r+3) & \cdots & (4r)(4r+1) \\
\cdots & \cdots & \cdots & \ddots & \cdots \\
2\cdots(2r+1) & 3\cdots(2r+2) & 4\cdots(2r+3) & \cdots &
(2r+2)\cdots(4r+1) \end{vmatrix} &
\end{matrix}.
\end{eqnarray*}
We obtain $A_{r}=\prod_{j=2}^{2r}j!$. Thus, there is a unique
solution for the system of nonhomogeneous linear equations:
\begin{eqnarray}
\left\{
   \begin{array}{ccccccccc}
   a_1 & + & a_2 & + & \cdots & + & a_{2r+1} & = &1, \\
   (2r+1)a_1 & + & (2r+2)a_2 & + & \cdots & + & (4r+1)a_{2r+1} & = &0, \\
   (2r+1)(2r)a_1 & + & (2r+1)(2r+2)a_2 & + & \cdots & + & (4r)(4r+1)a_{2r+1} & = &0, \\
   &&&\vdots&&&  \\
   (2r+1)!a_1 & + & 3 \cdots (2r+2)a_2 & + & \cdots & + & (2r+2) \cdots (4r+1)a_{2r+1} & =
   &0.
   \end{array} \right.\label{s2}
\end{eqnarray}
Let
\begin{eqnarray*}
\psi(x)=\left\{
\begin{array}{lrr}
a_1x^{2r+1}+a_2x^{2r+2}+\cdots+a_{2r+1}x^{4r+1}, &&0<x<1, \\
0,   &&x \leqslant0,  \\
1,  &&x=1.
             \end{array}
\right.
\end{eqnarray*}
with the coefficients $a_1,$ $a_2,$ $\cdots,$ $a_{2r+1}$ satisfying
(\ref{s2}). From (\ref{s2}), we see that $\psi(x)\in
C^{(2r)}(-\infty,+\infty)$, $0\leqslant\psi(x)\leqslant1$ for $0
\leqslant x\leqslant1$. Moreover, it holds that $\psi(1)=1$,\
$\psi^{(i)}(0)=0,\ i=0,1,\cdots,2r$. and $\psi^{(i)}(1)=0,\
i=1,2,\cdots,2r$. \\
Let $$H(f,x):=\sum_{i=1}^{r+1}f(x_{i})l_{i}(x),$$ and
$$l_{i}(x):=\frac{\prod_{j=1,j\neq i}^{r+1}(x-x_{j})}{\prod_{j=1,j\neq i}^{r+1}(x_{i}-x_{j})},\ x_{i}=\frac{[n\xi-({(r-1)/2}+i)]}{n},\ i=1,2, \cdots r+1.$$
Further, let
$$x_{1}^{'}=\frac{[n\xi-2\sqrt{n}]}{n},\ x_{2}^{'}=\frac{[n\xi-\sqrt{n}]}{n},\ x_{3}^{'}=\frac{[n\xi+\sqrt{n}]}{n},\ x_{4}^{'}=\frac{[n\xi+2\sqrt{n}]}{n},$$
and
$${\bar{\psi}}_{1}(x)=\psi(\frac{x-x_{1}^{'}}{x_{2}^{'}-x_{1}^{'}}),\ {\bar{\psi}}_{2}(x)=\psi(\frac{x-x_{3}^{'}}{x_{4}^{'}-x_{3}^{'}}).$$
Set
$${\bar{F}}_{n}(f,x):={\bar{F}}_{n}(x)=f(x)(1-{\bar{\psi}}_{1}(x)+{\bar{\psi}}_{2}(x))+{\bar{\psi}}_{1}(x)(1-{\bar{\psi}}_{2}(x))H(x).$$
We have
\begin{eqnarray*}
{\bar{F}}_{n}(f,x)=\left\{\begin{array}{lr}
f(x),          &       x\in [0,x_{r-5/2}]\cup [x_{r+3/2},1],   \\
f(x)(1-{\bar{\psi}}_{1}(x))+{\bar{\psi}}_{1}(x)
H(x),      &
x\in [x_{r-5/2},x_{r-3/2}],  \\
H(x),          &       x\in [x_{r-3/2},x_{r+1/2}],  \\
H(x)(1-{\bar{\psi}}_{2}(x))+{\bar{\psi}}_{2}(x)f(x), & x\in
[x_{r+1/2},x_{r+3/2}].
            \end{array}
\right.
\end{eqnarray*}
Obviously, ${\bar{F}}_{n}(f,x)$ is linear, reproduces polynomials of
degree $r$, and ${\bar{F}}_{n}(f,x)\in C^{(2r)}([0,1])$, provided
that $f \in C^{(2r)}([0,1]).$\\~\\
Now, we can define our new combinations of Bernstein operators as
follows:
$${\bar{B}}_{n,r}(f,x):=B_{n,r}({\bar{F}_{n}},x)=\sum_{i=0}^{r-1}C_{i}(n)B_{n_i}({\bar{F}_{n}},x),$$
where $C_{i}(n)$ satisfy the conditions (a)-(d). Our main result is
the following: \\~\\
\textbf{Theorem.} For any $\alpha >0,$ $0\leqslant \lambda
\leqslant1,$ we have
\begin{eqnarray}
\|{\bar{w}}{\bar{B}}_{n,r}^{(2r)}(f)\|\leqslant
Cn^{2r}\|{\bar{w}}f\|,\ f\in W_{\bar{w}}^{2r},\label{s3}
\end{eqnarray}
\begin{eqnarray}
|{\bar{w}(x)}\varphi^{2r\lambda}(x){\bar{B}}_{n,r}^{(2r)}(f,x)|\leqslant
\left\{
\begin{array}{lrr}
Cn^{r}\{\max\{n^{r(1-\lambda)},\varphi^{2r(\lambda-1)}\}\}\|{\bar{w}}f\|,    &&f\in C_{\bar{w}},    \\
C(\|{\bar{w}}f\|+\|{\bar{w}}\varphi^{2r\lambda}f^{(2r)}\|), &&f\in
W_{{\bar{w}}}^{2r},
              \end{array}
\right.\label{s4}
\end{eqnarray}
\begin{eqnarray}
\|{\bar{w}}{\bar{B}}_{n,r}(f)\|\leqslant C\|{\bar{w}}f\|,\ f\in
C_{\bar{w}},\label{s5}
\end{eqnarray}
\begin{eqnarray}
\|{\bar{w}}({\bar{B}}_{n,r}(f)-f)\|\leqslant \left\{
\begin{array}{lrr}
{\frac{C}{n^{r}}}
(\|{\bar{w}}f\|+\|{\bar{w}}\varphi^{2r}f^{(2r)}\|), &&f\in
W_{\bar{w}}^{2r},    \\
C(\omega_{\varphi}^{2r}(f,n^{-1/2})_{\bar{w}}+n^{-r}\|{\bar{w}}f\|),
&&f\in C_{\bar{w}},
              \end{array}
\right.\label{s6}
\end{eqnarray}
and for $0< \gamma <2r,$
\begin{eqnarray}
\|{\bar{w}}({\bar{B}}_{n,r}(f)-f)\|=O(n^{-\gamma/2})
\Longleftrightarrow
\omega_{\varphi}^{2r}(f,t)_{\bar{w}}=O(t^{r}).\label{s7}
\end{eqnarray}
\section{Lemmas}
\textbf{Lemma 1.}(\cite{Zhou}) For any non-negative real $u$ and
$v$, we have
\begin{eqnarray}
\sum_{k=1}^{n-1}({\frac kn})^{-u}(1-{\frac
kn})^{-v}p_{nk}(x)\leqslant Cx^{-u}(1-x)^{-v}.\label{s8}
\end{eqnarray}
\textbf{Lemma 2.}
For any positive real $\alpha $, and $f\in W_{\bar{w}}^{2r}$, we
have
\begin{eqnarray}
\|{\bar{w}}\varphi^{2r-2j}f^{(2r-j)}\|\leqslant
C(\|{\bar{w}}f\|+\|{\bar{w}}\varphi^{2r}f^{(2r)}\|).\label{s9}
\end{eqnarray}
\textbf{Proof.} \textit{Case 1.} $\xi \in[0,{\frac 14}]\cup [{\frac
34},1]$. It follows from Kolmogolov's inequality that
\begin{eqnarray*}
|f^{(2r-j)}({\frac 12})|\leqslant
C(\|f\|_{[1/4,3/4]}+\|f^{(2r)}\|_{[1/4,3/4]}),
\end{eqnarray*}
Moreover,
\begin{eqnarray}|f^{(2r-j)}({\frac 12})|\leqslant
C(\|{\bar{w}}f\|_{[1/4,3/4]}+\|{\bar{w}}\varphi^{2r}f^{(2r)}\|_{[1/4,3/4]}).\label{s10}
\end{eqnarray}
When $0\leqslant x\leqslant {\frac 12},$ $u$ between $x$ and
$\frac{k}{n},$ we have $\frac{|k/n-u|^{r-1}}{{\bar{w}}(u)}\leqslant
\frac{|k/n-x|^{r-1}}{{\bar{w}}(x)},$ then
\begin{eqnarray*}
|f^{(2r-j)}(x)-f^{(2r-j)}({\frac 12})|&\leqslant& \int_{x}^{\frac
12}|f^{(2r-j+1)}(u)|du\\
&\leqslant&C\|{\bar{w}}\varphi^{2r-2j+2}f^{(2r-j+1)}\|\int_{x}^{\frac
12}\frac{du}{{\bar{w}(u)}\varphi^{2r-2j+2}(u)}\\
&=&C\|{\bar{w}}\varphi^{2r-2j+2}f^{(2r-j+1)}\|\int_{x}^{\frac
12}\frac{|k/n-u|^{r-1}du}{|k/n-u|^{r-1}{\bar{w}(u)}\varphi^{2r-2j+2}(u)}\\
&\leqslant&C\|{\bar{w}}\varphi^{2r-2j+2}f^{(2r-j+1)}\|(\int_{x}^{\frac
12}\frac{|k/n-u|^{r-1}du}{|k/n-x|^{r-1}{\bar{w}(u)}\varphi^{2r-2j+2}(u)}\\&&+\int_{x}^{\frac
12}\frac{|k/n-u|^{r-1}du}{|k/n-1/2|^{r-1}{\bar{w}(u)}\varphi^{2r-2j+2}(u)})\\
&\leqslant&C\|{\bar{w}}\varphi^{2r-2j+2}f^{(2r-j+1)}\|\frac{x^{-r+j}}{{\bar{w}(x)}},
\end{eqnarray*}
which, together with (\ref{s10}), gives that
\begin{eqnarray*}
|{\bar{w}(x)}\varphi^{2r-2j}(x)f^{(2r-j)}(x)|\leqslant C(\|{\bar{w}}\varphi^{2r-2j+2}f^{(2r-j+1)}\|+
\|{\bar{w}}f\|+\|{\bar{w}}\varphi^{2r}f^{(2r)}\|).
\end{eqnarray*}
Similarly, we can prove that the above inequality also holds when
$1/2< x\leqslant1.$ Therefore, we obtain that
\begin{eqnarray}
|{\bar{w}(x)}\varphi^{2r-2j}(x)f^{(2r-j)}(x)|\leqslant
C(\|{\bar{w}}\varphi^{2r-2j+2}f^{(2r-j+1)}\|+
\|{\bar{w}}f\|+\|{\bar{w}}\varphi^{2r}f^{(2r)}\|).\label{s11}
\end{eqnarray}
Now, the result follows from (\ref{s11}) when $j=1$, and thus the
result can be deduced from (\ref{s11}) by induction when
$1<j\leqslant
r.$\\~\\
\textit{Case 2.} $\xi \in[{\frac 14},{\frac
34}]\cup\{\frac{1}{2}\}$. The situation goes similarly.\\~\\
\textbf{Lemma 3.} For any $f\in W_{\bar{w}}^{2r}$, we have
\begin{eqnarray}
\|{\bar{w}}(f-H)\|_{[x_{r-5/2},x_{r+3/2}]}\leqslant
{\frac{C}{n^{r}}}
(\|{\bar{w}}f\|+\|{\bar{w}}\varphi^{2r}f^{(2r)}\|).\label{s12}
\end{eqnarray}
\textbf{Proof.} By Taylor expansion, we have
\begin{eqnarray}
f({x_i})=\sum_{u=0}^{r}\frac{(x_i-x)^u}{u!}f^{(u)}(x)+{\frac
1{r!}}\int_{x}^{x_{i}}(x_i-s)^rf^{(r+1)}(s)ds,\label{s13}
\end{eqnarray}
It follows from (\ref{s13}) and the identity
\begin{eqnarray*}
\sum\limits_{i=1}^{r+1}x_{i}^{v}l_{i}(x)=x^v,\ v=0,1,\ldots,r.
\end{eqnarray*}
we have
\begin{eqnarray*}
H(f,x)&=&\sum_{i=0}^{r}\sum_{u=0}^{r}\frac{(x_i-x)^u}{u!}f^{(u)}(x)l_{i}(x)+{\frac
1{r!}}\sum_{i=1}^{r+1}l_{i}(x)\int_{x}^{x_{i}}(x_i-s)^rf^{(r+1)}(s)ds\nonumber\\
&=&f(x)+\sum_{u=1}^{r}f^{(u)}(x)(\sum_{v=0}^{u}C_{u}^{v}(-x)^{u-v}\sum_{i=1}^{r}x_{i}^{v}l_{i}(x))\nonumber\\
&&+{\frac
1{r!}}\sum_{i=1}^{r+1}l_{i}(x)\int_{x}^{x_{i}}(x_i-s)^rf^{(r+1)}(s)ds,
\end{eqnarray*}
which implies that
$${\bar{w}(x)}|f(x)-H(f,x)|={\frac 1{r!}}{\bar{w}(x)}\sum_{i=1}^{r+1}l_{i}(x)\int_{x}^{x_{i}}(x_i-s)^rf^{(r+1)}(s)ds,$$
since $|l_{i}(x)|\leqslant C$ for $x\in [0,{\frac
1n}],\ i=1,2,\cdots,r+1$.  \\
It follows from $\frac{|x_i-s|^{r-1}}{{\bar{w}}(s)}\leqslant
\frac{|x_i-x|^{r-1}}{{\bar{w}}(x)},$ $s$ between $x_i$ and $x$, then
\begin{eqnarray*}
{\bar{w}(x)}|f(x)-H(f,x)|&=&C\frac{\bar{w}(x)}{n^r}\sum_{i=1}^{r+1}\int_{x}^{x_{i}}|f^{(r+1)}(s)|ds\nonumber\\
&\leqslant&C\frac{\bar{w}(x)}{n^r}\|{\bar{w}}\varphi^{2}f^{(r+1)}\|\sum_{i=1}^{r+1}\int_{x}^{x_{i}}|{\bar{w}}^{-1}(s)\varphi^{-2}(s)ds|\nonumber\\
&\leqslant&{\frac C{n^r}}\|{\bar{w}}\varphi^{2}f^{(r+1)}\|.
\end{eqnarray*}
which, together with (\ref{s9}) (when $j=r-1$) implies (\ref{s12}).\\~\\
\textbf{Lemma 4.} For any $f\in W_{{\bar{w}}}^{2r}$ and $\alpha
>0$, we have
\begin{eqnarray}
\|{\bar{w}}\varphi^{2r}{\bar{F}_{n}^{(2r)}}\|\leqslant
C(\|{\bar{w}}\varphi^{2r}f^{(2r)}\|+ \|{\bar{w}}f\|).\label{s14}
\end{eqnarray}
\textbf{Proof.} We only prove the above result when $x\in
[x_{r-5/2},x_{r-3/2}]$, the others can be done similarly. Obviously,
\begin{eqnarray*}
|{\bar{F}_{n}^{(2r)}}(x)|&=&|(H(x)+{\bar{\psi}}_{1}(x)(f(x)-H(x))^{(2r)}|\nonumber\\
&=&|\sum_{i=0}^{2r}C_{2r}^{i}({\bar{\psi}}_{1}(x))^{(i)}(f(x)-H(x))^{(2r-i)}|\nonumber\\
&\leqslant& C\sum\limits_{i=0}^{2r}n^{\frac
i2}|(f(x)-H(x))^{(2r-i)}|
\end{eqnarray*}
If $2r-i\geqslant r+1$, using (\ref{s9}), then
\begin{eqnarray*}
|{\bar{w}(x)}\varphi^{2r}(x)(f(x)-H(x))^{(2r-i)}|&=&|{\bar{w}(x)}\varphi^{2r-2i}(x)f^{(2r-i)}(x)|\cdot
\varphi^{2i}(x)\nonumber\\
&\leqslant& Cn^{-i}(\|{\bar{w}}\varphi^{2r}f^{(2r)}\|+
\|{\bar{w}}f\|).
\end{eqnarray*}
If $0<2r-i<r+1$, by the following well-known inequality
\begin{eqnarray*}
\|g^{(j)}\|\leqslant
C((d-c)^{-j}\|g\|_{[c,d]}+(d-c)^{(2r-j)}\|g^{(2r)}\|_{[c,d]}),\
0<j<2r,
\end{eqnarray*}
we get
\begin{eqnarray*}
&&|{\bar{w}(x)}\varphi^{2r}(x)(f(x)-H(x))^{(2r-i)}|\nonumber\\&\leqslant&
C{\bar{w}(x)}\varphi^{2r}(x)(n^{2r-i}\|f-H\|_{[x_{r-5/2},x_{r-3/2}]}+n^{-i}\|f^{(2r)}\|_{[x_{r-5/2},x_{r-3/2}]})\nonumber\\
&\leqslant&
C(n^{r-i}\|{\bar{w}}(f-H)\|_{[x_{r-5/2},x_{r-3/2}]}+n^{-i}\|{\bar{w}}\varphi^{2r}f^{(2r)}\|_{[x_{r-5/2},x_{r-3/2}]})\nonumber\\
&\leqslant& Cn^{-i}(\|{\bar{w}}\varphi^{2r}f^{(2r)}\|+
_{}\|{\bar{w}}f\|).
\end{eqnarray*}
If $i=2r$, by (\ref{s12}), we have
\begin{eqnarray*}
|{\bar{w}(x)}\varphi^{2r}(x)(f(x)-H(x))^{(2r-i)}|&=&|{\bar{w}(x)}\varphi^{2r}(x)(f(x)-H(x))|\nonumber\\
&\leqslant& Cn^{-2r}(\|{\bar{w}}\varphi^{2r}f^{(2r)}\|+
\|{\bar{w}}f\|).
\end{eqnarray*}
Now the lemma follows from bringing these results together.\\~\\
\textbf{Lemma 5.} Let
$A_n(x):={\bar{w}(x)}\sum\limits_{|k-n\xi|\leqslant
\sqrt{n}}p_{n,k}(x)$. Then $A_n(x)\leqslant Cn^{-\alpha/2}$ for
$0<\xi <1$ and $\alpha>0$. \\~\\
\textbf{Proof.} If $|x- \xi|\leqslant {\frac {3}{\sqrt{n}} }$, then
the statement is trivial. Hence assume $0 \leqslant x \leqslant
\xi-{\frac {3}{\sqrt{n}} }$ (the case $\xi+{\frac {3}{\sqrt{n}} }
\leqslant x \leqslant 1$ can be treated similarly). Then for a fixed
$x$ the maximum of $p_{n,k}(x)$ is attained for
$k=k_n:=[n\xi-\sqrt{n}]$. By using Stirling's formula, we get
\begin{eqnarray*}
p_{n,k_n}(x)&\leqslant& C{\frac {({\frac
{n}{e}})^n\sqrt{n}x^{k_n}(1-x)^{n-k_n}}{({\frac
{k_n}{e}})^{k_n}\sqrt{k_n}({\frac
{n-k_n}{e}})^{n-k_n}\sqrt{n-k_n}}}\nonumber\\
&\leqslant& {\frac {C}{\sqrt{n}}}({\frac {nx}{k_n}})^{k_n}({\frac
{n(1-x)}{n-k_n}})^{n-k_n}\nonumber\\
&=&{\frac {C}{\sqrt{n}}}(1-{\frac {k_n-nx}{k_n}})^{k_n}(1+{\frac
{k_n-nx}{n-k_n}})^{n-k_n}.
\end{eqnarray*}
Now from the inequalities
$$k_n-nx=[n\xi-\sqrt{n}]-nx>n(\xi-x)-\sqrt{n}-1\geqslant {\frac 12}n(\xi-x),$$
and
$$1-u\leqslant e^{-u-{\frac 12}u^2},\ 1+u\leqslant e^u,\ u\geqslant 0.$$
We have that the second inequality is valid. To prove the first one
we consider the function $\lambda(u)=e^{-u-{\frac 12}u^2}+u-1.$ Here
$\lambda(0)=0,\ \lambda^\prime(u)=-(1+u)e^{-u-{\frac 12}u^2}+1,\
\lambda^\prime(0)=0,\ \lambda^{\prime\prime}(u)=u(u+2)e^{-u-{\frac
12}u^2}\geqslant 0,$ whence $\lambda(u)\geqslant 0$ for $u\geqslant
0$. Hence
\begin{eqnarray*}
p_{n,k_n}(x) &\leqslant&{\frac {C}{\sqrt{n}}}exp\{k_n[-{\frac
{k_n-nx}{k_n}}-{\frac 12}({\frac {k_n-nx}{k_n}})^2] +
k_n-nx\}\nonumber\\
&=&{\frac {C}{\sqrt{n}}}exp\{\frac {({k_n-nx})^2}{2k_n}\}\leqslant
e^{-Cn(\xi-x)^2}.
\end{eqnarray*}
Thus $A_n(x)\leqslant C(\xi-x)^\alpha e^{-Cn(\xi-x)^2}$. An easy
calculation shows that here the maximum is attained when
$\xi-x={\frac {C}{\sqrt{n}}}$ and the lemma follows.\\~\\
\textbf{Lemma 6.} For $0<\xi <1,\ \alpha, \beta>0$, we
have
\begin{eqnarray}
{\bar{w}(x)}\sum\limits_{|k-n\xi|\leqslant \sqrt{n}}|k-nx|^\beta
p_{n,k}(x)\leqslant
Cn^{(\beta-\alpha/2)}\varphi^\beta(x).\label{s15}
\end{eqnarray}
\textbf{Proof.} By (\ref{s8}) and the lemma 5, we have
\begin{eqnarray*}
{\bar{w}(x)}^{\frac
{1}{2n}}({\bar{w}(x)\sum\limits_{|k-n\xi|\leqslant
\sqrt{n}}p_{n,k}(x)})^{\frac
{2n-1}{2n}}(\sum\limits_{|k-n\xi|\leqslant \sqrt{n}}|k-nx|^{2n\beta}
p_{n,k}(x))^{\frac {1}{2n}}\leqslant
Cn^{(\beta-\alpha/2)}\varphi^\beta(x).
\end{eqnarray*}
\section{Proof of Theorem 1}
\subsection{Proof of (\ref{s3})}
We first prove $x\in [0,{\frac 1n}]$ (The same as $x\in [1-{\frac
1n},1]$), now
\begin{eqnarray}
|{\bar{w}}(x){\bar{B}}_{n,r}^{(2r)}(f,x)|&\leqslant&
{\bar{w}}(x)\sum_{i=0}^{r-1}{\frac
{n_{i}!}{({n_{i}-2r})!}}\sum_{k=0}^{n_i-2r}|C_{i}(n)\overrightarrow{\Delta}_{\frac
1{n_i}}^{2r}{\bar{F}}_{n}{(\frac k{n_i})}|p_{n_i-2r,k}(x)\nonumber\\
&\leqslant&
C{\bar{w}}(x)\sum_{i=0}^{r-1}n_{i}^{2r}\sum_{k=0}^{n_i-2r}|C_{i}(n)\overrightarrow{\Delta}_{\frac
1{n_i}}^{2r}{\bar{F}}_{n}{(\frac k{n_i})}|p_{n_i-2r,k}(x)\nonumber\\
&\leqslant&
C{\bar{w}}(x)\sum_{i=0}^{r-1}n_{i}^{2r}\sum_{k=0}^{n_i-2r}\sum_{j=0}^{2r}C_{2r}^{j}|C_{i}(n){\bar{F}}_{n}({\frac
{k+2r-j}{n_i}})|p_{n_i-2r,k}(x)\nonumber\\
&\leqslant&
C{\bar{w}}(x)\sum_{i=0}^{r-1}n_{i}^{2r}\sum_{j=0}^{2r}C_{2r}^{j}|C_{i}(n){\bar{F}}_{n}({\frac
{2r-j}{n_i}})|p_{n_i-2r,0}(x)\nonumber\\
&&+
C{\bar{w}}(x)\sum_{i=0}^{r-1}n_{i}^{2r}\sum_{j=0}^{2r}C_{2r}^{j}|C_{i}(n){\bar{F}}_{n}({\frac
{n_{i}-j}{n_i}})|p_{n_i-2r,n_i-2r}(x)\nonumber\\
&&+
C{\bar{w}}(x)\sum_{i=0}^{r-1}n_{i}^{2r}\sum_{k=1}^{n_i-2r-1}\sum_{j=0}^{2r}C_{2r}^{j}|C_{i}(n){\bar{F}}_{n}({\frac
{k+2r-j}{n_i}})|p_{n_i-2r,k}(x)\nonumber\\
&:=&H_1 +H_2 + H_3. \label{s16}
\end{eqnarray}
We have
\begin{eqnarray*}
H_1&\leqslant&
C{\bar{w}}(x)\sum_{i=0}^{r-1}n_{i}^{2r}(\sum_{j=0}^{2r-1}|C_{i}(n){\bar{F}}_{n}({\frac
{2r-j}{n_i}})|+|{\bar{F}}_{n}({0})|)p_{n_i-2r,0}(x)\nonumber\\
&\leqslant&
Cn^{2r}\|{\bar{w}}f\|\sum_{i=0}^{r-1}\sum_{j=0}^{2r-1}(\frac{n_{i}|x-\xi|}{2r-j-n_{i}\xi})^{\alpha}(1-x)^{{n_{i}}-2r}\nonumber\\
&\leqslant&
Cn^{2r}\|{\bar{w}}f\|\sum_{i=0}^{r-1}(n_{i}|x-\xi|)^{\alpha}(1-x)^{{n_{i}}-2r}\nonumber\\
&\leqslant& Cn^{2r}\|{\bar{w}}f\|.
\end{eqnarray*}
Similarly, we can get $H_2\leqslant Cn^{2r}\|{\bar{w}}f\|$ and
$H_3\leqslant
Cn^{2r}\|{\bar{w}}f\|$. \\~\\
When $x\in [{\frac 1n},1-{\frac 1n}],$ according to \cite{Totik}, we
have
\begin{eqnarray}
&&|{\bar{w}}(x){\bar{B}}_{n,r}^{(2r)}(f,x)|\nonumber\\
&=&|{\bar{w}}(x)B_{n,r}^{(2r)}({\bar{F}_{n}},x)|\nonumber\\
&=&{\bar{w}}(x)(\varphi^{2}(x))^{-2r}\sum_{i=0}^{r-1}\sum_{j=0}^{2r}|Q_{j}(x,n_i)C_{i}(n)|n_{i}^{j}\sum_{k/n_i\in
A}|(x-{\frac kn_{i}})^{j}{\bar{F}}_{n}({\frac kn_{i}})|p_{n_i,k}(x)\nonumber\\
&&+{\bar{w}}(x)(\varphi^{2}(x))^{-2r}\sum_{i=0}^{r-1}\sum_{j=0}^{2r}|Q_{j}(x,n_i)C_{i}(n)|n_{i}^{j}\sum_{x_2^{\prime}
\leqslant k/n_i\leqslant x_3^\prime}|(x-{\frac kn_{i}})^{j}H({\frac
kn_{i}})|p_{n_i,k}(x)\nonumber\\
&:=&\sigma_1+ \sigma_2.\label{s17}
\end{eqnarray}
Where $A:=[0,x_2^{\prime}]\cup [x_3^{\prime},1]$, $H$ is a linear
function. If ${\frac kn_{i}}\in A,$ when ${\frac
{\bar{w}(x)}{\bar{w}(\frac {k}{n_{i}})}}\leqslant C(1+n_i^{-{\frac
{\alpha}{2}}}|k-n_ix|^\alpha),$ we have $|k-n_{i}\xi|\geqslant
{\frac {\sqrt{n_{i}}}{2}}$, then
\\$Q_{j}(x,n_i)=(n_ix(1-x))^{[(2r-j)/2]},$ and $(\varphi^{2}(x))^{-2r}Q_{j}(x,n_i)n_{i}^{j}\leqslant
C(n_i/\varphi^{2}(x))^{r+j/2}.$ By (\ref{s15}), then
\begin{eqnarray*}
\sigma_1&\leqslant&
C{\bar{w}}(x)\sum_{i=0}^{r-1}\sum_{j=0}^{2r}|C_{i}(n)|(\frac
{n_{i}}{\varphi^{2}(x)})^{r+j/2}\sum_{k=0}^{n_{i}}|(x-{\frac
kn_{i}})^{j}{\bar{F}}_{n}({\frac kn_{i}})|p_{n_i,k}(x)\nonumber\\
&\leqslant&
C\|{\bar{w}}f\|\sum_{i=0}^{r-1}\sum_{j=0}^{2r}|C_{i}(n)|(\frac
{n_{i}}{\varphi^{2}(x)})^{r+j/2}\sum_{k=0}^{n_{i}}[1+n_i^{-{\frac
{\alpha}{2}}}|k-n_ix|^\alpha]|x-{\frac kn_{i}|^{j}}p_{n_i,k}(x)\nonumber\\
&:=&I_1 + I_2.
\end{eqnarray*}
By a simple calculation, we have $I_1\leqslant
Cn^{2r}\|{\bar{w}}f\|$. By ({\ref{s8}), then
$$I_2\leqslant C\|{\bar{w}}f\|\sum_{i=0}^{r-1}\sum_{j=0}^{2r}|C_{i}(n)|n_i^{-({{\frac {\alpha}{2}}}+j)}(\frac {n_{i}}
{\varphi^{2}(x)})^{j/2}\sum_{k=0}^{n_{i}}|k-n_ix|^{\alpha+j}p_{n_i,k}(x)\leqslant
Cn^{2r}\|{\bar{w}}f\|.$$ We note that $|H({\frac kn_i})|\leqslant
max(|H(x_1^\prime)|,|H(x_4^\prime)|):=H(a)$.\\~\\
If $x\in [x_1^\prime,x_4^\prime],$ we have ${\bar{w}(x)}\leqslant
{\bar{w}(a)}.$ So, if $x\in [x_1^\prime,x_4^\prime],$ then
$$\sigma_2\leqslant C{\bar{w}(a)}H(a)n^r\varphi^{-2r}(x)\leqslant Cn^{2r}\|{\bar{w}}f\|.$$
If $x\notin [x_1^\prime,x_4^\prime],$ then
${\bar{w}(a)}>n_i^{-{\frac {\alpha}{2}}},$ we have
$$\sigma_2\leqslant C{\bar{w}(a)}H(a)\varphi^{-2r}(x){\bar{w}}(x)\sum_{i=0}^{r-1}C_{i}(n)n_{i}^{r+{\frac {\alpha}{2}}}
\sum_{x_2^{\prime} \leqslant k/n_i\leqslant
x_3^\prime}p_{n_i,k}(x)\leqslant Cn^{2r}\|{\bar{w}}f\|.$$ It follows
from combining the above inequalities (\ref{s16}) and (\ref{s17})
that the theorem is proved.
\subsection{Proof of (\ref{s4})}
(1) When $f\in C_{\bar{w}}$, we discuss it as follows:\\~\\
\textit{Case 1.} If $0\leqslant \varphi(x)\leqslant {\frac
{1}{\sqrt{n}}}$, by $(\ref{s3})$, we have
\begin{eqnarray}
|{\bar{w}(x)}\varphi^{2r\lambda}(x){\bar{B}}_{n,r}^{(2r)}(f,x)|\leqslant
Cn^{-r\lambda}|{\bar{w}(x)}{\bar{B}}_{n,r}^{(2r)}(f,x)|\leqslant
Cn^{r(2-\lambda)}\|{\bar{w}}f\|.\label{s18}
\end{eqnarray} \textit{Case 2.} If
$\varphi(x)> {\frac {1}{\sqrt{n}}}$, we have
\begin{eqnarray*}
&&|{\bar{B}}_{n,r}^{(2r)}(f,x)|=|B_{n,r}^{(2r)}({\bar{F}_{n}},x)|\nonumber\\
&\leqslant&(\varphi^{2}(x))^{-2r}\sum_{i=0}^{r-1}\sum_{j=0}^{2r}|Q_{j}(x,n_i)C_{i}(n)|n_{i}^{j}\sum_{k=0}^{n_i}|(x-{\frac
kn_{i}})^{j}{\bar{F}}_{n}({\frac kn_{i}})|p_{n_i,k}(x),
\end{eqnarray*}
$Q_{j}(x,n_i)=(n_ix(1-x))^{[(2r-j)/2]},$ and
$(\varphi^{2}(x))^{-2r}Q_{j}(x,n_i)n_{i}^{j}\leqslant
C(n_i/\varphi^{2}(x))^{r+j/2}$. So
\begin{eqnarray}
&&|{\bar{w}(x)}\varphi^{2r\lambda}(x){\bar{B}}_{n,r}^{(2r)}(f,x)|\nonumber\\
&\leqslant&
C{\bar{w}(x)}\varphi^{2r\lambda}(x)\sum_{i=0}^{r-1}\sum_{j=0}^{2r}|C_{i}(n)|({\frac
{n_{i}}{\varphi^2(x)}})^{r+j/2}\sum_{k=0}^{n_i}|(x-{\frac
kn_{i}})^{j}{\bar{F}}_{n}({\frac kn_{i}})|p_{n_i,k}(x)\nonumber\\
&=&
C{\bar{w}(x)}\varphi^{2r\lambda}(x)\sum_{i=0}^{r-1}\sum_{j=0}^{2r}|C_{i}(n)|({\frac
{n_{i}}{\varphi^2(x)}})^{r+j/2}\sum_{k/n_i\in A}|(x-{\frac
kn_{i}})^{j}{\bar{F}}_{n}({\frac kn_{i}})|p_{n_i,k}(x)\nonumber\\
&&+
C{\bar{w}(x)}\varphi^{2r\lambda}(x)\sum_{i=0}^{r-1}\sum_{j=0}^{2r}|C_{i}(n)|({\frac
{n_{i}}{\varphi^2(x)}})^{r+j/2}\sum_{x_2^{\prime} \leqslant
k/n_i\leqslant x_3^\prime}|{(x-{\frac
kn_{i}})^{j}}{\bar{F}}_{n}({\frac
kn_{i}})|p_{n_i,k}(x)\nonumber\\
&:=&\sigma_1+ \sigma_2.\label{s19}
\end{eqnarray}
Where $A:=[0,x_2^{\prime}]\cup [x_3^{\prime},1]$. We can easily get
$\sigma_1\leqslant n^r\varphi^{2r(\lambda-1)}(x)\|{\bar{w}}f\|,$
$\sigma_2\leqslant n^r\varphi^{2r(\lambda-1)}(x)\|{\bar{w}}f\|.$
By bringing these facts (\ref{s18}) and (\ref{s19}) together, the theorem is proved.\\~\\
(2) When $f\in W_{{\bar{w}}}^{2r},$ we have
\begin{eqnarray}
B_{n,r}^{(2r)}({\bar{F}_{n}},x)=\sum_{i=0}^{r-1}C_{i}(n)n_{i}^{2r}\sum_{k=0}^{n_{i}-2r}\overrightarrow{\Delta}_{\frac
1{n_{i}}}^{2r}{\bar{F}}_{n}({\frac
kn_{i}})p_{n_i-2r,k}(x).\label{s20}
\end{eqnarray}
If $0<k<n_{i}-2r,$ we have
\begin{eqnarray}
|\overrightarrow{\Delta}_{\frac 1{n_{i}}}^{2r}{\bar{F}}_{n}({\frac
kn_{i}}) |\leqslant Cn_{i}^{-2r+1}\int_{0}^{\frac
{2r}{n_{i}}}|{\bar{F}}_{n}^{(2r)}({\frac kn_{i}}+u)|du,\label{s21}
\end{eqnarray}
If $k=0,$ we have
\begin{eqnarray}
|\overrightarrow{\Delta}_{\frac
1{n_{i}}}^{2r}{\bar{F}}_{n}(0)|\leqslant
Cn_{i}^{-r+1}\int_{0}^{\frac
{2r}{n_{i}}}u^{2r-1}|{\bar{F}}_{n}^{(2r)}(u)|du,\label{s22}
\end{eqnarray}
Similarly
\begin{eqnarray*}
|\overrightarrow{\Delta}_{\frac 1{n_{i}}}^{2r}{\bar{F}}_{n}({\frac
{n_{i}-2r}{n_{i}}})|&\leqslant& Cn_i^{-2r+1}\int_{1-{\frac
{2r}{n_{i}}}}^{1}(1-u)^{2r-1}|{\bar{F}}_{n}^{(2r)}(u)|du.
\end{eqnarray*}
By (\ref{s20}) and (\ref{s21}), we have
\begin{eqnarray*}
|{\bar{w}(x)}\varphi^{2r\lambda}(x){\bar{B}}_{n,r}^{(2r)}(f,x)|&\leqslant&
C{\bar{w}(x)}\varphi^{2r\lambda}(x)\sum_{i=0}^{r-1}|C_{i}(n)|n_{i}^{2r}\sum_{k=0}^{n_{i}-2r}|\overrightarrow{\Delta}_{\frac
1{n_{i}}}^{2r}{\bar{F}}_{n}({\frac kn_{i}})|p_{n_i-2r,k}(x)\nonumber\\
&=&C{\bar{w}(x)}\varphi^{2r\lambda}(x)\sum_{i=0}^{r-1}|C_{i}(n)|n_{i}^{2r}\sum_{k=1}^{n_{i}-2r-1}|\overrightarrow{\Delta}_{\frac
1{n_{i}}}^{2r}{\bar{F}}_{n}({\frac kn_{i}})|p_{n_i-2r,k}(x)\nonumber\\
&&+
C{\bar{w}(x)}\varphi^{2r\lambda}(x)\sum_{i=0}^{r-1}|C_{i}(n)|n_{i}^{2r}|\overrightarrow{\Delta}_{\frac
1{n_{i}}}^{2r}{\bar{F}}_{n}(0)|p_{n_i-2r,0}(x)\nonumber\\
&&+
C{\bar{w}(x)}\varphi^{2r\lambda}(x)\sum_{i=0}^{r-1}|C_{i}(n)|n_{i}^{2r}|\overrightarrow{\Delta}_{\frac
1{n_{i}}}^{2r}{\bar{F}}_{n}({\frac
{n_{i}-2r}{n_{i}}})|p_{n_i-2r,n_i-2r}(x)\nonumber\\
&:=&I_1 + I_2 + I_3.
\end{eqnarray*}
By (\ref{s21}), we have
\begin{eqnarray*}
I_1&\leqslant&
C{\bar{w}(x)}\varphi^{2r\lambda}(x)\sum_{i=0}^{r-1}|C_{i}(n)|n_{i}\sum_{k=1}^{n_{i}-2r-1}\int_{0}^{\frac
{2r}{n_{i}}}|{\bar{F}}_{n}^{(2r)}({\frac
kn_{i}}+u)|dup_{n_i-2r,k}(x)\nonumber\\
&=&C{\bar{w}(x)}\varphi^{2r\lambda}(x)\sum_{i=0}^{r-1}|C_{i}(n)|n_{i}\sum_{k/n_i\in
A}\int_{0}^{\frac {2r}{n_{i}}}|{\bar{F}}_{n}^{(2r)}({\frac
kn_{i}}+u)|dup_{n_i-2r,k}(x)\nonumber\\
&&+
C{\bar{w}(x)}\varphi^{2r\lambda}(x)\sum_{i=0}^{r-1}|C_{i}(n)|n_{i}\sum_{x_2^{\prime}
\leqslant k/n_i\leqslant x_3^\prime}\int_{0}^{\frac
{2r}{n_{i}}}|H_{n}^{(2r)}({\frac
kn_{i}}+u)|dup_{n_i-2r,k}(x)\nonumber\\
&:=&T_1 + T_2.
\end{eqnarray*}
Where $A:=[0,x_2^{\prime}]\cup [x_3^{\prime},1]$, $H$ is a linear
function. If ${\frac kn_{i}}\in A$, when ${\frac
{\bar{w}(x)}{\bar{w}(\frac {k}{n_{i}})}}\leqslant C(1+n_i^{-{\frac
{\alpha}{2}}}|k-n_ix|^\alpha),$ we have $|k-n_{i}\xi|\geqslant
{\frac {\sqrt{n_{i}}}{2}},$ by (\ref{s8}) and (\ref{s14}), then
\begin{eqnarray*}
T_1 &\leqslant&
C\|{\bar{w}}\varphi^{2r\lambda}F^{(2r)}\|{\bar{w}(x)}\varphi^{2r\lambda}(x)\sum_{i=0}^{r-1}|C_{i}(n)|n_{i}\sum_{k/n_i\in
A}\int_{0}^{\frac {2r}{n_{i}}}{\bar{w}}^{-1}({\frac
{k}{n_{i}}}+u)\varphi^{-2r\lambda}({\frac
{k}{n_{i}}}+u)dup_{n_i-2r,k}(x)\nonumber\\
&\leqslant&
C\|{\bar{w}}\varphi^{2r\lambda}F^{(2r)}\|\varphi^{2r\lambda}(x)\sum_{i=0}^{r-1}|C_{i}(n)|n_{i}\sum_{k=0}^{n_{i}}\int_{0}^{\frac
{2r}{n_{i}}}[1+n_i^{-{\frac
{\alpha}{2}}}|k-n_ix|^\alpha]\varphi^{-2r\lambda}({\frac
{k}{n_{i}}})dup_{n_i-2r,k}(x)
\nonumber\\
&\leqslant&
C\|{\bar{w}}\varphi^{2r\lambda}{\bar{F}_{n}^{(2r)}}\|\leqslant
C(\|{\bar{w}}f\|+\|{\bar{w}}\varphi^{2r\lambda}f^{(2r)}\|).
\end{eqnarray*}
Similarly, we can get $T_2\leqslant
C(\|{\bar{w}}f\|+\|{\bar{w}}\varphi^{2r\lambda}f^{(2r)}\|)$. So
$I_1\leqslant
C(\|{\bar{w}}f\|+\|{\bar{w}}\varphi^{2r\lambda}f^{(2r)}\|)$ and by
(\ref{s22}), we have
\begin{eqnarray*}
I_2&\leqslant&
C{\bar{w}(x)}\varphi^{2r\lambda}(x)\sum_{i=0}^{r-1}|C_{i}(n)|n_{i}^{2r}|\overrightarrow{\Delta}_{\frac
1{n_{i}}}^{2r}{\bar{F}}_{n}(0)|p_{n_i-2r,0}(x)\nonumber\\
&\leqslant&
C{\bar{w}(x)}\varphi^{2r\lambda}(x)\sum_{i=0}^{r-1}|C_{i}(n)|n_{i}^{r+1}\int_{0}^{\frac
{2r}{n_{i}}}u^{2r-1}|{\bar{F}}_{n}^{(2r)}(u)|dup_{n_i-2r,0}(x)\nonumber\\
&\leqslant&
C\|{\bar{w}}\varphi^{2r\lambda}{\bar{F}_{n}^{(2r)}}\|\sum_{i=0}^{r-1}(n_{i}x)^{r(1+\lambda)}(1-x)^{r\lambda}\leqslant
C\|{\bar{w}}\varphi^{2r\lambda}{\bar{F}_{n}^{(2r)}}\|\nonumber\\
&\leqslant&
C(\|{\bar{w}}f\|+\|{\bar{w}}\varphi^{2r\lambda}f^{(2r)}\|).
\end{eqnarray*}
Analogously, $I_3\leqslant
C(\|{\bar{w}}f\|+\|{\bar{w}}\varphi^{2r\lambda}f^{(2r)}\|)$, then
the theorem is proved.\\~\\
\textbf{Corollary 1.} If $\alpha>0$ and $\lambda=0$, we have
\begin{eqnarray*}
|{\bar{w}(x)}{\bar{B}}_{n,r}^{(2r)}(f,x)|\leqslant \left\{
\begin{array}{lrr}
Cn^{2r}\|{\bar{w}}f\|,    &&f\in C_{\bar{w}},    \\
C(\|{\bar{w}}f\|+\|{\bar{w}}f^{(2r)}\|), &&f\in W_{{\bar{w}}}^{2r}.
              \end{array}
\right.\label{eight}
\end{eqnarray*}
\textbf{Corollary 2.} If $\alpha>0$ and $\lambda=1$, we have
\begin{eqnarray*}
|{\bar{w}(x)}\varphi^{2r}(x){\bar{B}}_{n,r}^{(2r)}(f,x)|\leqslant
\left\{
\begin{array}{lrr}
Cn^{r}\|{\bar{w}}f\|,    &&f\in C_{\bar{w}},   \\
C(\|{\bar{w}}f\|+\|{\bar{w}}\varphi^{2r}f^{(2r)}\|), &&f\in
W_{{\bar{w}}}^{2r}.
              \end{array}
\right.\label{eight}
\end{eqnarray*}
\subsection{Proof of (\ref{s5})}
\begin{eqnarray*}
|{\bar{w}(x)}{\bar{B}}_{n,r}(f,x)|&=&|{\bar{w}(x)}B_{n,r}({\bar{F}_{n}},x)|\leqslant
{\bar{w}(x)}\sum_{i=0}^{r-1}\sum_{k=1}^{n_i-1}|C_{i}(n){\bar{F}}_{n}{(\frac
k{n_i})}|p_{n_i,k}(x)\nonumber\\
&&+
{\bar{w}(x)}\sum_{i=0}^{r-1}|C_{i}(n){\bar{F}}_{n}{(0)}|p_{n_i,0}(x)+
{\bar{w}(x)}\sum_{i=0}^{r-1}|C_{i}(n){\bar{F}}_{n}{(1)}|p_{n_i,n_i}(x)\nonumber\\
&:=&I_1 + I_2 + I_3.
\end{eqnarray*}
Analogously, the theorem can be proved easily.
\subsection*{4.4. Proof of (\ref{s6})}
We assume $f\in W_{\bar{w}}^{2r},$ then
$\|{\bar{w}}({\bar{B}}_{n,r}(f)-{\bar{F}_{n}})\|\leqslant
{\frac{C}{n^{r}}}
(\|{\bar{w}}f\|+\|{\bar{w}}\varphi^{2r}f^{(2r)}\|).$
\\Recall that \cite{Totik}, then
\begin{eqnarray}
B_{n,r}((t-x)^j,x)=0,\ j=1,2,\cdots,r,\label{s23}
\end{eqnarray}
\begin{eqnarray}
B_{n,r}((t-x)^{2r-j},x)=O(n^{-r}\varphi^{2r-2j}(x)),\ x\in [{\frac
1n},1-{\frac 1n}],\ j=0,1,2,\cdots,r.\label{s24}
\end{eqnarray}
\textit{Case 1.} $x\in [{\frac 1n},1-{\frac 1n}].$ By using Taylor
expansion, we have
\begin{eqnarray*}
&&{\bar{w}(x)}({\bar{F}_{n}(x)}-B_{n,r}({\bar{F}_{n}},x))\nonumber\\
&=&{\bar{w}(x)}\sum_{j=1}^{2r-1}{\frac
{1}{(2r-j)!}}B_{n,r}((t-x)^{2r-j},x){\bar{F}_{n}^{(2r-j)}(x)}\nonumber\\
&&+ {\bar{w}(x)}B_{n,r}({\frac
{1}{(2r-j)!}}\int_x^t(t-u)^{2r-1}{\bar{F}_{n}^{(2r)}(u)}du,x)\nonumber\\
&:=&I_1 + I_2.
\end{eqnarray*}
By (\ref{s9}), (\ref{s14}) and (\ref{s24}), we have for $1\leqslant
j\leqslant r$, then
\begin{eqnarray}
{\frac
{{\bar{w}(x)\varphi^{2r-2j}(x)}}{n^r}}{\bar{F}_{n}^{(2r-j)}(x)}\leqslant
{\frac{C}{n^{r}}}
(\|{\bar{w}}{\bar{F}_{n}}\|+\|{\bar{w}}\varphi^{2r}{\bar{F}_{n}}^{(2r)}\|)\leqslant
{\frac{C}{n^{r}}}
(\|{\bar{w}}f\|+\|{\bar{w}}\varphi^{2r}f^{(2r)}\|),\label{s25}
\end{eqnarray}
By (\ref{s23}) and (\ref{s25}), we have
$$I_1\leqslant {\bar{w}(x)}\sum_{j=1}^{r-1}{\frac {1}{(2r-j)!}}|B_{n,r}((t-x)^{2r-j},x){\bar{F}_{n}^{(2r-j)}(x)}|\leqslant {\frac{C}{n^{r}}}
(\|{\bar{w}}f\|+\|{\bar{w}}\varphi^{2r}f^{(2r)}\|).$$ If $u$ is
between $t$ and $x$, we have
$\frac{|t-u|^{2r-1}}{\varphi^{2r}(u)}\leqslant
\frac{|t-x|^{2r-1}}{\varphi^{2r}(x)}$. Then
\begin{eqnarray*}
&&|{\bar{w}(x)}B_{n,r}({\frac
{1}{(2r-j)!}}\int_x^t(t-u)^{2r-1}{\bar{F}_{n}^{(2r)}(u)}du,x)|\nonumber\\
&\leqslant&
C{\bar{w}(x)}\sum_{i=0}^{r-1}\sum_{k=0}^{n_i}|C_{i}(n)|\int_x^{\frac
{k}{n_i}}|({\frac
{k}{n_i}}-u)^{2r-1}{\bar{F}_{n}^{(2r)}(u)}|dup_{n_i,k}(x)\nonumber\\
&=&C{\bar{w}(x)}\sum_{i=0}^{r-1}\sum_{k=1}^{n_i-1}|C_{i}(n)|\int_x^{\frac
{k}{n_i}}|({\frac
{k}{n_i}}-u)^{2r-1}{\bar{F}_{n}^{(2r)}(u)}|dup_{n_i,k}(x)\nonumber\\
&&+
C{\bar{w}(x)}\sum_{i=0}^{r-1}|C_{i}(n)|(1-x)^{n_i}\int_0^x{u^{2r-1}}|{\bar{F}_{n}^{(2r)}(u)}|du\nonumber\\
&&+
C{\bar{w}(x)}\sum_{i=0}^{r-1}|C_{i}(n)|x^{n_i}\int_x^1{(1-u)^{2r-1}}|{\bar{F}_{n}^{(2r)}(u)}|du\nonumber\\
&:=&J_1 + J_2 + J_3.
\end{eqnarray*}
We have
\begin{eqnarray*}
J_1&\leqslant&
C{\bar{w}(x)}\varphi^{-2r}(x)\sum_{i=0}^{r-1}\sum_{k/n_i\in
A}|C_{i}(n)({\frac {k}{n_i}}-x)^{2r-1}|\int_x^{\frac
{k}{n_i}}\varphi^{2r}(v)|{\bar{F}_{n}^{(2r)}(v)}|dvp_{n_i,k}(x)\nonumber\\
&&+ C{\bar{w}(x)}\varphi^{-2r}(x)\sum_{i=0}^{r-1}\sum_{x_2^{\prime}
\leqslant k/n_i\leqslant x_3^\prime}|C_{i}(n)({\frac
{k}{n_i}}-x)^{2r-1}|\int_x^{\frac
{k}{n_i}}\varphi^{2r}(v)|H^{(2r)}(v)|dvp_{n_i,k}(x)\nonumber\\
&:=&\sigma_1 + \sigma_2.
\end{eqnarray*}
Analogously, we can get $\sigma_1\leqslant {\frac{C}{n^{r}}}
(\|{\bar{w}}f\|+\|{\bar{w}}\varphi^{2r}f^{(2r)}\|)$. We note that
$|\varphi^{2r}(v)H^{(2r)}(v)|\leqslant
\max(|\varphi^{2r}(x_1^\prime)H^{(2r)}(x_1^\prime)|,\
|\varphi^{2r}(x_4^\prime)H^{(2r)}(x_4^\prime)|):=|\varphi^{2r}(a)H^{(2r)}(a)|,
$
$H^{(2r)}(x)$ is a linear function.\\~\\
If $x\in [x_1^\prime,x_4^\prime],$ then ${\bar{w}(x)}\leqslant
{\bar{w}(a)}.$ So, we have
\begin{eqnarray*}
\sigma_2&\leqslant&
C{\bar{w}(a)}\varphi^{2r}(a)|H^{(2r)}(a)|\varphi^{-2r}(x)\sum_{i=0}^{r-1}\sum_{k=1}^{n_i-1}|C_{i}(n)|({\frac
{k}{n_i}}-x)^{2r}p_{n_i,k}(x)\nonumber\\
&\leqslant& {\frac{C}{n^{r}}}
(\|{\bar{w}}f\|+\|{\bar{w}}\varphi^{2r}f^{(2r)}\|),
\end{eqnarray*}
If $x\notin [x_1^\prime,x_4^\prime]$, by ${\bar{w}(a)}>n_i^{-{\frac
{\alpha}2}}$, we have
\begin{eqnarray*}
\sigma_2&\leqslant&
C{\bar{w}(a)}\varphi^{-2r}(a)|H^{(2r)}(a)|\sum_{i=0}^{r-1}\sum_{x_2^{\prime}
\leqslant k/n_i\leqslant x_3^\prime}n_i^{\frac
{\alpha}2}|C_{i}(n)|({\frac {k}{n_i}}-x)^{2r}p_{n_i,k}(x)\nonumber\\
&\leqslant& {\frac{C}{n^{r}}}
(\|{\bar{w}}f\|+\|{\bar{w}}\varphi^{2r}f^{(2r)}\|).
\end{eqnarray*}
For $J_2,$ we have
\begin{eqnarray*}
J_2&\leqslant&
C\|{\bar{w}}\varphi^{2r}{\bar{F}_{n}}^{(2r)}\|{\bar{w}(x)}
\sum_{i=0}^{r-1}|C_{i}(n)|(1-x)^{n_i}\int_0^x{u^{2r-1}}{\bar{w}^{-1}(u)}\varphi^{-2r}(u)du\nonumber\\
&\leqslant& {\frac{C}{n^{r}}}
(\|{\bar{w}}f\|+\|{\bar{w}}\varphi^{2r}f^{(2r)}\|).
\end{eqnarray*}
\\Similarly, we have
$$J_3 \leqslant {\frac{C}{n^{r}}}
(\|{\bar{w}}f\|+\|{\bar{w}}\varphi^{2r}f^{(2r)}\|).$$
\\By bringing these facts together, we have
$$\|{\bar{w}}({\bar{B}}_{n,r}(f)-{\bar{F}_{n}})\|\leqslant
{\frac{C}{n^{r}}}
(\|{\bar{w}}f\|+\|{\bar{w}}\varphi^{2r}f^{(2r)}\|).$$
\textit{Case 2.} $x\in [0,{\frac 1n}]$ (Similarly as $x\in [1-{\frac
1n},1]$). By using Taylor expansion, we have
\begin{eqnarray*}
{\bar{w}(x)}|B_{n,r}({\bar{F}_{n}},x)-{\bar{F}_{n}(x)}|&\leqslant&
{\frac {\bar{w}(x)}{r!}}\sum_{i=0}^{r-1}|C_{i}(n)|B_{n_i}
(\int_x^t|(t-u)^r{\bar{F}_{n}}^{(r+1)}(u)|du,x)\nonumber\\
&&+ {\frac {\bar{w}(x)}{r!}}\sum_{i=0}^{r-1}|C_{i}(n)|(1-x)^{n_i}\int_0^xu^{2r-1}|{\bar{F}_{n}}^{(r+1)}(u)|du\nonumber\\
&:=&J_1 + J_2.
\end{eqnarray*}
\begin{eqnarray*}
J_1&\leqslant&
C{\bar{w}(x)}\sum_{i=0}^{r-1}\sum_{k=0}^{n_i}\int_x^{\frac
{k}{n_i}}|C_{i}(n)({\frac
{k}{n_i}}-u)^r{\bar{F}_{n}}^{(r+1)}(u)|dup_{n_i,k}(x)\nonumber\\
&:=&C{\bar{w}(x)}\sum_{i=0}^{r-1}\sum_{k=1}^{n_i-1}\int_x^{\frac {k}{n_i}}|C_{i}(n)({\frac {k}{n_i}}-u)^r{\bar{F}_{n}}^{(r+1)}(u)|dup_{n_i,k}(x)\nonumber\\
&&+ C{\bar{w}(x)}\sum_{i=0}^{r-1}|C_{i}(n)|x^{n_i}\int_x^{1}(1-u)^r|{\bar{F}_{n}}^{(r+1)}(u)|du\nonumber\\
&&+ C{\bar{w}(x)}\sum_{i=0}^{r-1}|C_{i}(n)|(1-x)^{n_i}\int_0^{x}u^r|{\bar{F}_{n}}^{(r+1)}(u)|du\nonumber\\
&:=&I_1 + I_2 + I_3.
\end{eqnarray*}
Analogously, we can get
\begin{eqnarray*}
I_1&\leqslant& {\frac{C}{n^{r}}}
(\|{\bar{w}}f\|+\|{\bar{w}}\varphi^{2r}f^{(2r)}\|),\nonumber\\
I_2&\leqslant& {\frac{C}{n^{r}}}
(\|{\bar{w}}f\|+\|{\bar{w}}\varphi^{2r}f^{(2r)}\|),\nonumber\\
I_3&\leqslant& {\frac{C}{n^{r}}}
(\|{\bar{w}}f\|+\|{\bar{w}}\varphi^{2r}f^{(2r)}\|).
\end{eqnarray*}
\begin{eqnarray}
J_1&\leqslant& {\frac{C}{n^{r}}}
(\|{\bar{w}}f\|+\|{\bar{w}}\varphi^{2r}f^{(2r)}\|),\nonumber\\
J_2&\leqslant& {\frac{C}{n^{r}}}
(\|{\bar{w}}f\|+\|{\bar{w}}\varphi^{2r}f^{(2r)}\|).\nonumber\\
\end{eqnarray}
So, we have
\begin{eqnarray*}
\|{\bar{w}}({\bar{B}}_{n,r}(f)-{\bar{F}_{n}})\|\leqslant
{\frac{C}{n^{r}}}
(\|{\bar{w}}f\|+\|{\bar{w}}\varphi^{2r}f^{(2r)}\|).
\end{eqnarray*}
Then
\begin{eqnarray*}
\|{\bar{w}}({\bar{B}}_{n,r}(f)-f)\|&\leqslant&
\|{\bar{w}}(f-{\bar{F}_{n}}(f))\|+\|{\bar{w}}({\bar{F}_{n}}(f)-{\bar{B}}_{n,r}(f))\|\nonumber\\
&\leqslant& {\frac{C}{n^{r}}}
(\|{\bar{w}}f\|+\|{\bar{w}}\varphi^{2r}f^{(2r)}\|).
\end{eqnarray*}
If $f\in C_{\bar{w}},$ there exists $g\in W_{\bar{w}}^{2r},$ by
(\ref{s5}) and the first inequality of (\ref{s6}), we have
\begin{eqnarray*}
\|{\bar{w}}({\bar{B}}_{n,r}(f)-f)\|&\leqslant& \|{\bar{w}}(f-g)\| +
\|{\bar{w}}{\bar{B}}_{n,r}(f-g)\| +
\|{\bar{w}}(g-{\bar{B}}_{n,r}(g))\|\nonumber\\
&\leqslant& C(\|{\bar{w}}(f-g)\|+
{\frac{1}{n^{r}}}(\|{\bar{w}}g\|+\|{\bar{w}}\varphi^{2r}g^{(2r)}\|))\nonumber\\
&\leqslant&
C(\omega_{\varphi}^{2r}(f,n^{-1/2})_{\bar{w}}+n^{-r}\|{\bar{w}}f\|).
\end{eqnarray*}
\subsection*{4.5. Proof of (\ref{s7})}
From the proof of (\ref{s6}), we actually have
$$\|{\bar{w}}({\bar{B}}_{n,r}(f)-f)\|\leqslant CK_{2r,\varphi}(f,t^r)_{\bar{w}}.$$
Therefore, $K_{2r,\varphi}(f,n^{-r})_{\bar{w}}=O({t^\alpha})$
implies
$$\|{\bar{w}}({\bar{B}}_{n,r}(f)-f)\|\leqslant (n^{-\alpha/2}).$$
By (\ref{s4}) and (\ref{s5}), we may choose g properly such that
$\|{\bar{w}}\varphi^{2r}g^{(2r)}\|<\infty$ and
\begin{eqnarray*}
\omega_{\varphi}^{2r}(f,n^{-1/2})_{\bar{w}} + {\frac
{\|{\bar{w}}f\|}{n^r}}&\leqslant&
\|{\bar{w}}({\bar{B}}_{n,r}(f)-f)\|+{\frac
1{n^r}}({\|\bar{w}}\varphi^{2r}{\bar{B}}_{n,r}^{(2r)}(f-g)\|
\nonumber\\
&&+\|{\bar{w}}\varphi^{2r}{\bar{B}}_{n,r}^{(2r)}(g)\|)+{\frac
{\|{\bar{w}}f\|}{n^r}}\nonumber\\
&\leqslant& \|{\bar{w}}(f-{\bar{B}}_{n,r}(f))\|+{\frac
{\|{\bar{w}}f\|}{n^r}}+C({\frac
kn})^r(\|{\bar{w}}(f-g)\|\nonumber\\
&&+k^{-r}{\|\bar{w}}\varphi^{2r}g^{(2r)}\|+k^{-r}\|{\bar{w}}f\|)\nonumber\\
&\leqslant& \|{\bar{w}}(f-{\bar{B}}_{n,r}(f))\|+C({\frac
kn})^r(\omega_{\varphi}^{2r}(f,k^{-1/2})_{\bar{w}}\nonumber\\
&&+k^{-r}\|{\bar{w}}f\|).
\end{eqnarray*}
Hence, by \cite{Totik}, we obtain the converse inequality.

\end{document}